\documentclass[12pt]{amsart}
\usepackage{amsthm}

\newcommand{\bz}{\mathbf Z}

\newcommand{\bc}{\mathbf C}

\newcommand{\br}{\mathbf R} 
\newcommand{\bh}{\mathbf H} 

\newcommand{\ind}{\operatorname{ind}}

\newcommand{\spin}{\text{spin}}

\newcommand{\Int}{\operatorname{Int}}
\newcommand{\coker}{\operatorname{coker}}

\newcommand{\Pin}{\operatorname{Pin}} 
\newcommand{\supp}{\operatorname{supp}}

\newtheorem{theorem}{Theorem}
\newtheorem{proposition}{Proposition}
\newtheorem{corollary}{Corollary}
\newtheorem{lemma}{Lemma}
\theoremstyle{definition}
\newtheorem{definition}{Definition}
\newtheorem{remark}{Remark}

\newtheorem{example}{Example}

\begin{document}
\title[ ]
{Constraints on intersection forms of spin 4-manifolds bounded by Seifert rational homology 3-spheres in terms of $\overline{\mu}$ and $\kappa$ invariants. 
}

\author[]{Masaaki Ue}

\address{
Kyoto University, Kyoto, 606-8502, Japan}
\email{ue.masaaki.87m@st.kyoto-u.ac.jp}

\begin{abstract}
We give some constraints on intersection forms of spin 4-manifolds bounded by Seifert rational homology 3-spheres in terms of 
$\overline{\mu}$-invariant, and compare them with those in terms of $\kappa$-invariant given by Manolescu. 
Furthermore in case of a Seifert rational homology 3-sphere, we show that 
the difference between $\kappa$ and $-\overline{\mu}$ is at most 2 (and there is an example where the difference is 2), 
and they coincide with each other under some extra conditions, 
including Floer $K_G$ split cases.
\end{abstract}

\maketitle

Let $Y$ be a closed 3-manifold  with spin structure $\mathfrak s$.  
Then $(Y, \mathfrak s )$ bounds a
compact smooth spin 4-manifold $(W,\mathfrak s_W )$, where $s_W$ is a spin 
structure of $W$ which coincides with $\mathfrak s$ on $Y$.
In the case where $Y$ is a Seifert (or more generally a plumbed) 
rational homology 3-sphere, there is a lift $\overline{\mu} (Y, \mathfrak s)$
of the Rokhlin invariant $\mu (Y, \mathfrak s)$, 
which is called the $\overline{\mu}$-invariant defined by 
Neumann (\cite{N})  and Siebenmann  (\cite{Si}).

In this paper we give the following constraints on the intersection forms of compact spin 4-manifolds bounded by Seifert rational homology 3-spheres 
in terms of the $\overline{\mu}$ invariant. In this paper we only consider Seifert fibrations over orientable base 2-orbifolds. 

\begin{theorem}\label{main1}
Let $(W, \mathfrak s_W )$ be a compact spin 4-manifold bounded by a Seifert rational homology 3-sphere with $\spin$ structure $(Y, \mathfrak s )$. 
Then  $\overline{\mu} (Y, \mathfrak s ) \equiv \sigma (W) /8 \pmod {2\bz}$ and 
\[
-b_2^- (W) +\frac{\sigma (W)}{8} \le \overline{\mu} (Y, \mathfrak s ) \le b_2^+ (W) +\frac{\sigma (W)}{8} 
\]
where $b_2^+ (W)$ (resp.~$b_2^- (W)$) is the dimension of the maximal positive definite (resp.~negative definite) subspace of the intersection form of 
$W$, and $\sigma (W)$ is the signature of $W$. 
\end{theorem}

There are similar constraints on the intersection form of $W$ in terms of the $\kappa$-invariant $\kappa (Y, \mathfrak s )$ 
 of $(Y, \mathfrak s )$ proved by Manolescu (which are valid for general rational homology 3-spheres). In a certain case $\kappa (Y, \mathfrak s )$ coincides with 
 $-\overline{\mu} (Y, \mathfrak s )$, but they are not necessarily the same even in the case of Brieskorn homology 3-spheres (\cite{Mano4}). 
 However by comparing the constraints in terms of $\kappa (Y, \mathfrak s )$ with the estimates given by $\overline{\mu} (Y, \mathfrak s )$, 
 we obtain the following result. 
 If $Y$ is an integral homology 3-sphere, we denote these invariants simply by $\kappa (Y)$ and $\overline{\mu} (Y)$ since
 $Y$ has a unique $\spin$ structure.

 \begin{theorem}\label{main2}
 Let $(Y, \mathfrak s )$ be a Seifert rational homology 3-sphere over a base 2-orbifold of genus 0 with $\spin$ structure. 
 
 \begin{enumerate}
 \item
 We have the following inequality. 
 \[
 0 \le \kappa (Y , \mathfrak s) +\overline{\mu} (Y, \mathfrak s ) \le 2 .
 \]
 \item 
 If $(Y, \mathfrak s )$ is Floer $K_G$ split, then $\kappa (Y, \mathfrak s ) +\overline{\mu} (Y, \mathfrak s ) =0$.  
\item 
 We have $0 \le \kappa (Y, \mathfrak s) +\kappa (-Y , \mathfrak s ) \le 4$ (although $\kappa (Y, \mathfrak s ) +\kappa (-Y , \mathfrak s )$ is not 
 necessarily $0$). 

 \item If one of the multiplicities of the singular fibers for the Seifert fibration of $Y$ is even and $\deg Y >0$, then 
 $\kappa (Y, \mathfrak s ) =-\overline{\mu} (Y, \mathfrak s)$. 
 
  \item 
 If $Y$ is a Seifert integral homology 3-sphere (with the unique $\spin$ structure), 
 $\kappa (Y) +\overline{\mu} (Y) =0$ or $2$. The same claim holds for $-Y$, and 
 hence 
 $\kappa (Y ) +\kappa (-Y) =0$, $2$, or $4$. 
   \end{enumerate}
 \end{theorem}
 
 Theorem \ref{main2} (4) can be extended to the case when all multiplicities are odd under certain extra conditions on $(Y, \mathfrak s )$ (Remark \ref{addenda}). 
 
 \section{The Fukumoto-Furuta invariant and the proof of Theorem 1.}
 
 The proof of Theorem 1 is based on the estimates of the Fukumoto-Furuta invariant (\cite{FF}) and its relation to the $\overline{\mu}$-invariant of a Seifert 
 rational homology 3-sphere, which we recall in this section.
 
 \begin{definition}
 Let $(X, \mathfrak s_X )$ be a 
 compact $\spin$ 4-orbifold with $\spin$ structure $\mathfrak s_X $ bounded by
 a rational homology 3-sphere $(Y, \mathfrak s )$ with $\spin$ structure. 
 Choose a compact spin 4-manifold $(X' , \mathfrak s_{X'} )$ 
 with $\partial (X' , \mathfrak s_{X'}  ) =(Y, \mathfrak s )$ and put $(Z, \mathfrak s_Z ) = (X\cup (-X') , \mathfrak s_X \cup 
 \mathfrak s_{X '} )$. Then the Fukumoto-Furuta invariant $w(Y, X, \mathfrak s_X )$  is defined to be 
 \[
 w(Y, X, \mathfrak s_X ) =-\ind \mathcal D_Z  (\mathfrak s_Z ) +\frac{\sigma (X' )}{8} 
 \]
 where $\mathcal D_Z (\mathfrak s_Z )$ is the Dirac operator of $(Z, \mathfrak s_Z )$. 
 The value of $w(Y, X, \mathfrak s_X )$ does not depend on the choice of $X'$. 
\end{definition}

We note that our sign convention of $w(Y, X, \mathfrak s_X )$ is opposite to the original definition by Fukumoto-Furuta \cite{FF} and \cite{FFU}. 
 Let $\mathrm{Sing} X$ be the set of all singularities of $X$. Hereafter we only consider 4-orbifolds with isolated singularities. In such cases 
 the link of every $x\in \mathrm{Sing}  X$ is a spherical 3-manifold with spin structure induced from $\mathfrak s_X$. 
 Then the regular neighborhood of $x_i \in \mathrm{Sing}  X$ in $X$ is a cone $cS_i$ over  a spherical 3-manifold $S_i$, and 
$X_0 =X \setminus \cup_{x_i \in {\mathrm Sing} X} \Int cS_i$
 is a 4-manifold with $\spin$ structure $\mathfrak s_{X_0}$ induced from $\mathfrak s_X$ satisfying 
 $\partial X_0 =Y \cup \cup_{i} (-S_i )$. 
 
 \begin{proposition}\label{orbifold index} \cite{FFU}
The orbifold index theorem shows that 
 \begin{align*}
 \ind \mathcal D_Z (\mathfrak s_Z ) &=-\frac{1}{24} \int_Z p_1 (Z) +\sum_{i} \delta^{\mathrm{Dir}} (S_i , \mathfrak s_{S_i} ), \\
 \sigma (Z) &=\frac 13 \int_Z p_1 (Z) +\sum_i \delta^\mathrm{sign} (S_i ) 
 \end{align*}
where $p_1 (Z)$ is the first Pontrjagin class of $Z$, $\mathfrak s_{S_i}$ is the $\spin$ structure of $S_i$ induced from $\mathfrak s_X$, 
$\delta^{\mathrm{Dir}} (S_i  , \mathfrak s_{S_i} )$ and $\delta^{\mathrm{sign}} (S_i )$ are the contributions of $x_i \in  \mathrm{Sing} X$ to 
the corresponding indices. Using these formulae $w (Y, X, \mathfrak s_X )$ is represented as follows. 
\[
w(Y, X, \mathfrak s_X ) =\frac 18 (\sigma (Z) +\sum_{x_i \in \mathrm{Sing} X} \delta_{x_i}  ) +\frac{\sigma (X' )}{8} =
\frac 18 (\sigma (X) +\sum_{x_i \in \mathrm{Sing} X } \delta_{x_i}). 
\]
Here we put $\delta_{x_i} =-(8\delta^{\mathrm{Dir}} (S_i , \mathfrak s_{S_i} ) +\delta^{\mathrm{sign}} (S_i ))$. 

 \end{proposition}
 
 For later use, we reformulate the contributions of $\mathrm{Sing} X$ to the above indices in terms of the $\eta$ invariants (\cite{APS}). 
 
 \begin{proposition}\label{U3} \cite{U3}
 Let $g_i$ be a standard metric of $S_i$ (which is induced from the metric of the constant curvature on the 3-sphere of radius 1 as the universal covering of 
 $S_i$). Then we have 
 \[
 \eta^{\mathrm{Dir}} (S_i , \mathfrak s_{S_i} , g_i ) =2\delta^{\mathrm{Dir}} (S_i , \mathfrak s_{S_i} ) , \quad 
 \eta^{\mathrm{sign}} (S_i , g_i ) =\delta^{\mathrm{sign}} (S_i ) .
 \]
 Here $\eta^{\mathrm{Dir}} (S_i , \mathfrak s_{S_i} , g_i )$ and $\eta^{\mathrm{sign}} (S_i , g_i )$ are the $\eta$ invariants of 
 the Dirac operator of $(S_i , \mathfrak s_{S_i} , g_i )$ and the sigunature operator of $(S_i , g_i )$ respectively. 
 We note that the kernel of the Dirac operator of $(S_i , \mathfrak s_{S_i} , g_i )$ is $0$ since $g_i$ is a  metric of positive scalar curvature. 
 \end{proposition} 
 
 The value of $w(Y, X, \mathfrak s_X )$ depends on both $Y$ and $X$ in general, but if $Y$ is a Seifert rational homology 3-sphere, it
 is related to the $\overline{\mu}$ invariant of $Y$ by some particular choice of $X$. 
 
 \begin{definition}
 Let $P(\Gamma )$ be a plumbed 4-manifold associated with a weighted graph $\Gamma$.  We assume that a plumbed 3-manifold $Y =\partial P(\Gamma )$ 
 is a rational homology 3-sphere. Thus $\Gamma$ is a tree such that each vertex $v$ with weight $e(v)$ of $\Gamma$ corresponds to a 
 $D^2$-bundle over $S^2$ with Euler class $e(v)$. 
 Let $\{ e_v \}$ be a basis of  $H_2 (P(\Gamma ), \bz )$ each of which is  represented by a zero section of the $D^2$-bundle over $S^2$ corresponding to $v$. 
If we fix a $\spin$ structure $\mathfrak s_Y$ of $Y=\partial P(\Gamma )$, 
 there exists a unique characteristic element $w_{\Gamma} =\sum_v \epsilon_v e_v $ satisfying the following conditions. 
 \begin{itemize}
 \item $\epsilon_v$ is either $0$ or $1$. 
 \item The Poincar\'e dual $PD w_{\Gamma} \pmod 2 \in H^2 (P(\Gamma ) , \partial P(\Gamma ), \bz_2 )$ of 
$w_{\Gamma}$ is a lift of $w_2 (P (\Gamma ))$, which is an obstruction to 
 extending $\mathfrak s_Y$ to a $\spin$ structure of $P(\Gamma )$. 
 \end{itemize}
 
 We denote by $\supp ( \Gamma)$ the set $\{ v \in \Gamma \ | \ \epsilon_v =1 \}$. 
 Then the $\overline{\mu}$ invariant of $(Y, \mathfrak s_Y )$ is defined to be 
 \[
 \overline{\mu} (Y, \mathfrak s_Y ) =\frac 18 (\sigma (P(\Gamma )) -w_{\Gamma} \cdot w_{\Gamma} ) .
 \]
 The value of $\overline{\mu} (Y, \mathfrak s_Y )$ does not depend on $\Gamma$, although $\Gamma$ satisfying 
 $Y= \partial P(\Gamma )$ is not unique. 
 \end{definition}
 
 \begin{remark}\label{other invariants}
 The $\overline{\mu}$ invariant has the following properties. 
 \begin{enumerate}
 \item $\overline{\mu} (Y, \mathfrak s_Y )\equiv \mu (Y, \mathfrak s_Y ) \mod{2\bz}$. If $Y$ is a plumbed integral homology 3-sphere with 
 the unique $\spin$ structure, 
 $\overline{\mu} (Y)$ is an integral lift of $\mu (Y) \in \bz_2$. 
 \item If $Y$ is a Seifert rational homology 3-sphere (with base 2-orbifold of genus $0$), then $Y$ is a plumbed 3-manifold corresponding to a 
 star-shaped graph. 
 \item 
 In some cases the $\overline{\mu}$ invariant is related to other $\spin$ cobordism invariants. 
For example, the relations between $\overline{\mu} (Y, \mathfrak s )$ and the correction term $d(Y, \mathfrak s )$ of the Heegaard Floer homology (\cite{OS}) 
are discussed in \cite{U3}, \cite{Stip}.  If $Y$ is a Seifert integral homology 3-sphere, $-\overline{\mu} (Y)$ coincides with the $\beta$-invariant $\beta (Y)$
defined by Manolescu \cite{Mano3} (Stoffregen \cite{Stoff}), and 
the same equality holds  if $Y$ is a plumbed  rational homology 3-sphere corresponding to a graph with at most one bad vertex (Dai \cite{Dai}). 
In case of a Seifert integral homology 3-sphere, the relations between $\overline{\mu} (Y )$ and the invariants $\alpha (Y)$, $\gamma (Y)$ and  $d(Y)$ are 
also determined in \cite{Stoff}. 
More generally, $\underline{d} (Y, \mathfrak s ) = -2\overline{\mu} (Y, \mathfrak s )$ for every almost rational plumbed rational homology 3-sphere $Y$
(Dai-Manolescu \cite{DM}), where $\underline{d} (Y, \mathfrak s )$ is the correction term of the involutive Heegaard Floer homology (\cite{HM}). 
In all above cases, the $\overline{\mu}$ invariant is a $\spin$ rational homology cobordism invariant. 
 \end{enumerate}
 \end{remark}
 
 The relation between the $\overline{\mu}$ invariant and the Fukumoto-Furuta invariant for a Seifert fibration is discussed in \cite{FF}, \cite{FFU}, \cite{S}, \cite{U1}, \cite{U2}.

 If a Seifert rational homology 3-sphere $Y$ has a fibration over a 2-orbifold $S^2 (a_1 , \cdots , a_n )$ of genus $0$ with singular points of multiplicity 
 $a_i$ $(1 \le i \le n )$, then $Y$ is represented by Seifert invariants of the form 
$
 \{ b, (a_1 , b_1 ) , \dots , (a_n , b_n ) \}
$
 satisfying $0< b_i < a_i$, $\gcd (a_i , b_i ) = 1$ $(1\le i \le n )$, and 
 the degree of $Y$ is defined to be $\deg Y = b+\sum_{i=1}^n b_i /a_i$ (we follow the definition of \cite{MOY}, whose sign convention is opposite to 
 that in \cite{U2}, where we follow 
 the definition of \cite{NR}). 
 In case of a Seifert rational homology 3-sphere, we have the following result. 
 
 \begin{proposition}{\label{U2}} \cite{U2}{\footnote
 {The definitions of $w$ and $\overline{\mu}$ invariants in \cite{U2} and  \cite{U3} are both 8 times those in this paper.}
 }
 Let $(Y, \mathfrak s )$ be a Seifert rational homology 3-sphere  with $\spin$ structure. Then there exist compact
 spin 4-orbifolds $(X^{\pm} , \mathfrak s_{X^{\pm}} )$ 
 with $\partial (X^{\pm} , \mathfrak s_{X^{\pm}} ) =(Y, \mathfrak s )$ 
 satisfying the following conditions. 
 \begin{enumerate}
 \item $b_1 (X^{\pm} ) =0$, $b_2^+ (X^+ ) \le  1$, and $b_2^- (X^- ) \le 1$. 
 \item 
$ w(Y, X^+ , \mathfrak s_{X^+ } ) =w(Y, X^- , \mathfrak s_{X^- } ) =\overline{\mu} (Y, \mathfrak s )$. 
\end{enumerate}
If $Y$ is spherical, then we can choose $X^{\pm}$ so that both of them are the cone over $Y$. 
Suppose that one of the multiplicities of the singular fibers for the Seifert fibration of $Y$ is even. Then 
we can choose $X^{\pm}$ so that $X^+ =X^-$, satisfying $b_2^+ (X^{\pm} ) =1$ and $b_2^- (X^{\pm}) =0$ if $\deg Y >0$, and 
$b_2^+ (X^{\pm}  )=0$ and $b_2^- (X^{\pm} ) =1$ if $\deg Y <0$. 
\end{proposition}
 
 To construct $X^{\pm}$, we choose appropriate weighted graphs $\Gamma^{\pm}$ satisfying $\partial P(\Gamma^{\pm} ) =Y$ and subgraphs 
 $\Gamma^{\pm}_0$ of $\Gamma^{\pm}$ such that $\Gamma^{\pm}_0$ contain $\supp (\Gamma^{\pm} )$ and 
  each of $S^{\pm} := \partial P (\Gamma_0^{\pm} )$ is a union of spherical 3 -manifolds (in fact lens spaces). 
 Then embed $P(\Gamma^{\pm}_0 )$ in $\Int P(\Gamma^{\pm} )$ and 
 replace their images with the cones over $S^{\pm}$ to obtain $X^{\pm}$. 
 Such $X^{\pm}$ are constructed by Saveliev \cite{S} for Seifert integral homology 3-spheres, and the result is extended to the case of 
 Seifert rational homology 3-spheres in \cite{U2}. 
 
 The proof of Theorem \ref{main1} is based on Proposition \ref{U2} and the following orbifold $10/8$-theorem. 
 
 \begin{theorem}{\label{8/10}} \cite{FF}
 Let $(Z, \mathfrak s_Z )$ be a closed spin 4-orbifold with $\spin$ structure $\mathfrak s_Z$. Then 
 \begin{enumerate}
 \item $\ind \mathcal D_Z (\mathfrak s_Z ) \equiv 0 \pmod 2$, 
 \item $\ind \mathcal D_Z (\mathfrak s_Z ) =0$ or otherwise 
 \[
 1-b_2^- (Z) \le \ind \mathcal D_Z (\mathfrak s_Z ) \le b_2^+ (Z) -1 .
 \]
 \end{enumerate}
 \end{theorem}
 
 For a compact spin manifold $(W, \mathfrak s_W )$ with $\partial (W, \mathfrak s_W ) =(Y, \mathfrak s )$, let 
 $(Z^{\pm}, \mathfrak s_{Z^{\pm}} ) = (X^{\pm} , \mathfrak s_{X^{\pm}} ) \cup (-W , \mathfrak s_W )$ be closed $\spin$ 4-orbifolds where 
 $X^{\pm}$ are the 4-orbifolds chosen in Proposition \ref{U2}. 
 Then applying Theorem \ref{8/10} to $(Z^{\pm}, \mathfrak s_{Z^{\pm}})$ we have $\overline{\mu} (Y, \mathfrak s ) \equiv \frac{\sigma (W)}{8} 
 \pmod {2\bz}$ and 
 
 \begin{align*}
 \overline{\mu} (Y, \mathfrak s) &=w( Y , X^- , \mathfrak s_{X^-} ) =-\ind \mathcal D_{Z^-} (\mathfrak s_{Z^-} ) +\frac{\sigma (W)}{8} \le 
b_2^- (Z^- ) -1 +\frac{\sigma (W)}{8}  \\
&=b_2^- (X^-  ) +b_2^+ (W) -1 +\frac{\sigma (W)}{8} \le b_2^+ (W) +\frac{\sigma (W)}{8},  \\
\overline{\mu} (Y, \mathfrak s) &=w(Y, X^+ , \mathfrak s_{X^+} ) =-\ind \mathcal D_{Z^+} (\mathfrak s_{Z^+} ) +
\frac{\sigma (W)}{8} \ge 1- b_2^+ (Z^+ ) +\frac{\sigma (W)}{8} \\
&=1-b_2^- (W) -b_2^+ (X^+ ) +\frac{\sigma (W)}{8} 
 \ge - b_2^- (W) +\frac{\sigma (W)}{8}. 
 \end{align*}
 
 The above inequalities hold even if $\ind \mathcal D_{Z^{\pm} } (\mathfrak s_{Z^{\pm}} ) =0$. 
 The second inequality is also deduced from the first one with respect to $X^{\pm}$ chosen for $-Y$ and $-W$, since $\overline{\mu} (-Y ) =-\overline{\mu} (Y)$. 
 This completes the proof of Theorem \ref{main1}. 
We will use Proposition \ref{U2} to prove Theorem \ref{main2} in \S 3. 
 
 \begin{remark}\label{remark on main1}
 In some cases we obtain slightly better estimates than those in Theorem \ref{main1} by the above inequalities. If $Y$ is a spherical 3-manifold $S$, we can choose 
 $X^{\pm}$ so that $X^+ =X^- =cS$ and hence applying  Theorem \ref{8/10} to $Z =cS \cup (-W)$ we have either 
 \[
 1-b_2^- (W) +\frac{\sigma (W)}{8} \le \overline{\mu} (S, \mathfrak s ) \le b_2^+ (W) -1 +\frac{\sigma (W)}{8} 
 \]
and $\overline{\mu} (S, \mathfrak s ) \equiv \frac{\sigma(W)}{8} \pmod{2\bz}$, or $\overline{\mu} (S, \mathfrak s ) =\frac{\sigma (W)}{8}$.

If $Y$ has a fibration over $S^2 (a_1 , \dots , a_n )$ and one of $a_i$ is even, then 
we can choose $X^{\pm}$ so that $X^+ =X^-$ as in Proposition \ref{U2}. If $\deg Y >0$, $X^{\pm}$ satisfies $b_2^+ (X^{\pm} ) =1$, 
$b_2^- (X^{\pm} ) =0$, and hence  we have either 
\[
-b_2^- (W) +\frac{\sigma (W)}{8} \le \overline{\mu} (Y, \mathfrak s ) \le b_2^+ (W) -1 +\frac{\sigma (W)}{8} 
\]
or $\overline{\mu} (Y, \mathfrak s ) =\frac{\sigma(W)} {8}$. If $\deg Y <0$, then since $X^{\pm}$ satisfies $b_2^+ (X^{\pm} ) =0$, 
$b_2^- (X^{\pm} ) =1$, we have either 
 \[
 1-b_2^- (W) +\frac{\sigma  (W)}{8} \le \overline{\mu} (Y, \mathfrak s ) \le b_2^+ (W) +\frac {\sigma (W)}{8}
 \]
or $\overline{\mu} (Y, \mathfrak s ) =\frac{\sigma (W)}{8}$. In either case we have $\overline{\mu} (Y, \mathfrak s ) \equiv 
\frac{\sigma(W)}{8} \pmod{2\bz}$. 
 If $Y$ is a Seifert integral homology 3-sphere, then 
 $Y$ is uniquely determined by $a_i$'s up to orientation, and is denoted by $\Sigma (a_1 , \dots , a_n )$ if $\deg Y<0$. 
 Then $\Sigma (a_1 , \dots , a_n )$ bounds an associated orbifold disk bundle $D (a_1 , \dots , a_n )$ over $S^2 (a_1 , \dots , a_n )$. 
 If one of $a_i$ is even, then $D (a_1 , \dots , a_n )$ is $\spin$ and $b_2^- (D (a_1 , \dots , a_n ) ) =1$, $b_2^+ (D (a_1 , \dots , a_n ) =0$
 (\cite{FF}). 
 Hence putting $X^{\pm } =D^2 (a_1 , \dots , a_n )$ we also obtain the same inequality. 
But in other cases $D (a_1 , \dots , a_n )$ is not necessarily spin, so we cannot use it.  
  \end{remark}

\section{The $\kappa$ invariant for a rational homology 3-sphere.}

For a rational homology 3-sphere, $K$-theoretic invariants are defined via the theory of 
the $\Pin (2)$ Seiberg-Witten-Floer stable homotopy type (\cite{Mano1}, \cite{Mano2}) in \cite{FL}, \cite{Mano4}. 
In this section we recall the definition of Manolescu's $\kappa$ invariant for a rational homology 3-sphere with $\spin$ structure $(Y, \mathfrak s )$ and 
its properties \cite{Mano4}
{\footnote{Theorems \ref{estimate1}, \ref{estimate2}, and Corollary \ref{cor} below are described for 
integral homology 3-spheres in \cite{Mano4}. But as is pointed out in \cite{Mano4}, these results are easily generalized to the case of 
rational homology 3-spheres. We give the outline of their proofs below.}}
Hereafter we put $G=\Pin (2) =S^1 \cup j S^1$. 

\begin{definition} \cite{Mano4} 
\begin{enumerate}
 \item Let $\widetilde{\br}$ (resp.~ $\widetilde{\bc}$) be a real (resp.~complex) 1-dimensional $G$- representation space where 
 $G$ acts as  multiplication by $\pm 1$ via the natural projection $G\to G/S^1 \cong \{ \pm 1 \}$, and let $\bh$ be the  set of quaternions on which $G$ 
 acts as left multiplication. Then the representation ring $R (G)$ is generated by $z=2-\bh$ and $w =1-\widetilde{\bc}$ with relations 
 $w ^2 =wz =2w$. 
 \item
 A compact $G$-space $X$ is called a space of type SWF at level $s$ if 
 the fixed point set $X^{S^1}$ of the action of $S^1 \subset G$ on $X$ is $G$ homotopy equivalent to 
 $( \widetilde{\br}^s )^+$ and $G$ acts freely on $X \setminus X^{S^1}$. 
 \end{enumerate}
 \end{definition}
 
 \begin{definition} \cite{Mano4}
 Let $X$ be a space of type SWF at level $2t$ (in which case $X^{S^1}$ is $G$ homotopy equivalent to $(\widetilde{\bc}^t )^+$). 
 Then the inclusion $\iota : X^{S^1} \to X$ induces the map $\iota^* : \widetilde K_G (X) \to \widetilde K_G (X^{S^1} )$ whose image is 
 represented as $\mathcal J (X) \cdot \beta_{ (\widetilde{\bc}^t )^+ }$, where $\mathcal J (X)$ is an ideal of 
 $R(G) \cong \widetilde K_G (S^0 )$ and $\beta_{(\widetilde{\bc}^t )^+ }$ is the Bott class of $\widetilde K_G ( (\widetilde{\bc}^t )^+ )$. 
 
 \begin{enumerate}
 \item
 There exists $k\ge 0$ such that both $w^k $and $z^k$ belong to $\mathcal J (X)$, so we can define 
 \[
 k(X) =\min \{ k\ge 0 \ | \ \text{there exists $x \in \mathcal J (X)$ such that $wx=z^k w $} \} 
 \]
 \item $X$ is called $K_G$ split if $\mathcal J  (X)=(z^k )$ for some $k\ge 0$. 
 \end{enumerate}
 \end{definition}
 
\begin{lemma}\label{estimate for SWF} \cite{Mano4} 
\begin{enumerate}
\item
Let $f: X\to X'$ is a $G$ equivariant map between two spaces of type SWF at level $2t$, and assume that $f^{S^1} : X^{S^1} \to {X'}^{S^1}$ is a 
$G$ homotopy equivalence, where $f^{S^1}$ is the restriction of $f$ to $X^{S^1}$. Then $k(X) \le k (X')$. 
\item Let $X$ and $X'$ be spaces of type SWF at level $2t$ and $2t '$ respectively and $t<t'$. Let $f: X\to X'$ be a $G$ equivariant map such that 
$f^G : X^G \to {X'}^G$ is a homotopy equivalence, where $f^G$ is the restriction of $f$ to the fixed point set of the $G$ action on $X$. 
Then $k(X) +t \le k(X') +t'$. Furthermore, if $X$ is $K_G$ split, then $k(X) +t+1 \le k(X' ) +t'$. 
\end{enumerate}
\end{lemma}
 
 A space of type SWF associated with a rational homology 3-sphere with $\spin$ structure $(Y, \mathfrak s )$ is constructed as follows. We fix a Riemann metric 
 $g$ of $Y$. Let $V = i \ker d^* \oplus \Gamma (S )$ where $i\ker d^* \subset i \Omega^1 (Y)$ and $\Gamma (S)$ is the space of sections of 
 the Spinor bundle $S$ over $Y$ associated with $\mathfrak s$. 
 Then the action of $G$ on $V$ is defined as 
 \[
 e^{i\theta} (a, \phi ) =(a, e^{i\theta} \phi ) , \quad j (a, \phi ) = (-a , j\phi ) , 
 \]
and the gradient of the Chern-Simons-Dirac functional $\nabla CSD$ over $V$ is represented as a sum of a $G$ equivariant linear operator $\ell$ and 
a compact operator $c$. Let $V^{\mu}_{\tau}$ be the finite dimensional subspace of $V$ spanned by eigenvectors corresponding to 
the eigenvalues of $\ell$ belonging to $(\tau , \mu ]$. Then $V^{\mu}_{\tau}$ is a direct sum of 
$V^{\mu}_{\tau} (\widetilde{\br} ) =\widetilde{\br}^s$ and $V^{\mu}_{\tau} (\bh ) =\bh^t$ for some $s$, $t$.  
Manolescu defined a finite dimensional approximation of $\nabla CSD$ over $V^{\nu}_{-\nu}$ and corresponding flows 
(approximated Seiberg-Witten flows) 
for a large $\nu >0$ and showed that 
a  set $S_{\nu}$ of critical points and flows between them in a large ball centered at the origin in $V^{\nu}_{-\nu}$ forms an isolated invariant set. 
Thus a $G$ equivariant Conley index $I_{\nu}$ of $S_{\nu}$ is defined up to $G$ homotopy equivalence (\cite{Mano1}, \cite{Mano3}, \cite{Mano4}). 

\begin{proposition} \cite{Mano3}, \cite{Mano4} 
\begin{enumerate}
\item There exists a unique critical point $\Theta$ of the above flow corresponding to the unique isolated and nondegenerate 
reducible critical point of 
$\nabla CSD$ (which is perturbed if necessary). The Conley index $I(\Theta )$ of $\Theta$ is $(V^0_{-\nu} )^+$.
\item $I_{\nu}$ is a space of type SWF at level $\dim V^0_{-\nu} (\widetilde{\br} )$. 
\end{enumerate}

\end{proposition}

\begin{definition} \cite{Mano3}
For $(Y, \mathfrak s )$ choose a compact 4-manifold with $\spin$ structure $(W , \mathfrak s_W )$ bounded by 
$(Y, \mathfrak s )$ and a metric $g_W$ of $W$ extending the given metric $g$ of $Y$. Then define $n(Y, \mathfrak s , g )$ to be 
\[
n(Y, \mathfrak s , g ) =\ind_{\bc} \mathcal D_W (\mathfrak s_W ) +\frac{\sigma (W)}{8} \in \frac 18 \bz  
\]
where $\ind_{\bc} \mathcal D_W (\mathfrak s_W )$ is the Atiyah-Patodi-Singer index of the Dirac operator of $(W, \mathfrak s_W )$ 
associated with $g_W$. 
The value of $n(Y, \mathfrak s, g )$ is independent of the choice of $(W, \mathfrak s_W)$. If $Y$ is an integral homology 3-sphere, 
$n(Y, \mathfrak s , g ) \in \bz$ since $\sigma (W)$ is divisible by 8. 
\end{definition}

\begin{definition} \cite{Mano1} \cite{Mano3} \cite{Mano4} 
A Seiberg-Witten Floer stable homotopy type 
$SWF (Y, \mathfrak s )$ of $(Y, \mathfrak s )$ is a formal desuspension of $I_{\nu}$ (up to stable equivalence) defined as follows:
If $\dim V^0_{-\nu} (\widetilde{\br} )$ is even, then $V^0_{-\nu} =\widetilde{\br}^{2k} \oplus \bh^{\ell}$ for some $k$, $\ell$ and 
\[
SWF (Y, \mathfrak s ) = \Sigma^{-(k\widetilde{\bc} +(\ell +\frac 12 n(Y, \mathfrak s , g ) )\bh )} I_{\nu}, 
\]
while if 
$\dim V^0_{\-nu} (\widetilde{\br} )$ is odd, then 
$V^0_{-\nu} =\widetilde{\br}^{2k-1} \oplus \bh^{\ell}$ for some $k$, $\ell$ and 
\[
SWF (Y, \mathfrak s ) =\Sigma^{-(k\widetilde{\bc} +(\ell +\frac 12 n(Y, \mathfrak s , g ) )\bh )} (\Sigma^{\widetilde{\br}} I_{\nu} ) .
\]

More precisely, $SWF (Y, \mathfrak s )$ is an object 
of the $G$ equivariant graded suspension category $\mathfrak C$.  
\end{definition}

\begin{definition} \cite{Mano4}
The $\kappa$ invariant of $(Y, \mathfrak s )$ is defined to be 
\begin{align*}
\kappa (Y, \mathfrak s ) &=2 k (SWF (Y, \mathfrak s ) )  \\
&= 
\begin{cases} 
2k(I_{\nu} ) -2\dim_{\bh} V^0_{-\nu} (\bh ) -n(Y,\mathfrak s , g ) & \text{if $I_{\nu}$ is at even level, } \\
2k(\Sigma^{\widetilde{\br}} I_{\nu} ) -2\dim_{\bh} V^0_{-\nu} (\bh ) -n(Y, \mathfrak s , g ) & \text{if $I_{\nu}$ is at odd level}.
\end{cases} 
\end{align*}

\end{definition}

\begin{proposition}\label{homology cobordism} \cite{Mano4} 
\begin{enumerate}
\item 
$\kappa (Y, \mathfrak s ) \equiv \overline{\mu} (Y, \mathfrak s ) \pmod {2\bz}$. If $(Y, \mathfrak s )$ is an integral homology 3-sphere with the unique 
$\spin$ structure, then $\kappa (Y, \mathfrak s )$ is an integral lift of $\overline{\mu} (Y, \mathfrak s ) \in \bz_2$. 

\item 
$\kappa (Y, \mathfrak s )$ is a $\spin$ rational homology cobordism invariant. 
\end{enumerate}
\end{proposition}

For a $\spin$ 4-manifold $(W, \mathfrak s_W )$ bounded by $(Y, \mathfrak s )$,
$\ker \mathcal D_W (\mathfrak s_W )$ and $\coker \mathcal D_W (\mathfrak s_W)$ are vector spaces over $\bh$, and
$\ind_{\bc} \mathcal D_W (\mathfrak s_W )
=2\ind_{\bh} \mathcal D_W (\mathfrak s_W )$.
It follows that 
$\kappa (Y, \mathfrak s ) \equiv n(Y, \mathfrak s , g ) \equiv \sigma (W) /8 \equiv \mu (Y, \mathfrak s ) \pmod{2\bz}$. 
The second claim in Proposition \ref{homology cobordism} is proved by applying Theorem \ref{estimate1} below to a $\spin$ rational homology 
cobordism between rational homology 3-spheres.

\begin{theorem}\label{estimate1}\cite{Mano4} 

Let $(Y_i , \mathfrak s_i )$ $(i=0,1 )$ be rational homology 3-spheres with $\spin$ structures and 
$(W , \mathfrak s_W )$ be a $\spin$ cobordism with $b_1 (W) =0$ from $(Y_0 , \mathfrak s_0 )$ to $(Y_1 , \mathfrak s_1 )$. 
Then 
\[
-b_2^- (W) +\frac{\sigma (W)}{8} -1 \le \kappa (Y_0, \mathfrak s_0 ) -\kappa (Y_1 , \mathfrak s_1 ) \le b_2^+ (W) +\frac{\sigma (W)}{8} +1 . 
\]
The first term of the above inequality can be replaced by $-b_2^- (W) +\frac{\sigma (W)}{8}$ if $b_2^- (W)$ is even, and 
the last term can be replaced by $b_2^+ (W) +\frac{\sigma (W)}{8}$ if $b_2^+ (W)$ is even. 
\end{theorem}

\begin{definition} \cite{Mano4} 
$(Y, \mathfrak s )$ is called Floer $K_G$ split if $I_{\nu}$ 
(resp.~$\Sigma^{\widetilde{\br}} I_{\nu}$) is $K_G$ split when $I_{\nu}$ is at even level (resp.~at odd level), where $I_{\nu}$ is the Conley index 
that appears in the definition of $\kappa (Y, \mathfrak s )$. 
\end{definition}

If one of $(Y_i , \mathfrak s_i )$ is Floer $K_G$ split, the inequalities in Theorem \ref{estimate1} are slightly improved. 

\begin{theorem}\label{estimate2} \cite{Mano4}
Let $(W, \mathfrak s_W )$ be a $\spin$ cobordism as in Theorem \ref{estimate1}.  
\begin{enumerate}
\item 
If $(Y_0 , \mathfrak s_0 )$ is Floer $K_G$ split, then 
\begin{align*}
&\kappa (Y_0, \mathfrak s_0 ) -\kappa  (Y_1 , \mathfrak s_1 )   \\
&\le \begin{cases} 
b_2^+ (W) +\frac{\sigma (W)}{8} -1 &\quad \text{if $b_2^+ (W)$ is odd}, \\
b_2^+ (W) +\frac{\sigma (W)}{8} -2 & \quad \text{if $b_2^+ (W)$ is even and $b_2^+ (W) >0$. } 
\end{cases}
\end{align*}

\item 
If $(Y_1 , \mathfrak s_1 )$ is Floer $K_G$ split, then 
\begin{align*}
&\kappa (Y_1 , \mathfrak s_1 ) -\kappa (Y_0 , \mathfrak s_0 )   \\
&\le \begin{cases}
b_2^- (W) -\frac{\sigma (W) }{8} -1 &\quad \text{if $b_2^- (W)$ is odd, } \\
b_2^- (W) -\frac{\sigma (W)} {8}  -2 &\quad \text{if $b_2^- (W)$ is even and $b_2^- (W) > 0$. }
\end{cases}
\end{align*}
\end{enumerate}
\end{theorem}

Since $S^3$ with the unique $\spin$ structure is Floer $K_G$ split and $\kappa (S^3 ) =0$, we have the following corollary. 

\begin{corollary}\label{cor} 
\cite{Mano4}

Let $(Y, \mathfrak s )$ be a rational homology 3-sphere with $\spin$ structure and $(W, \mathfrak s_W )$ be a compact $\spin$ 4-manifold 
with $b_1 (W)=0$ and $\partial (W, \mathfrak s_W ) =(Y, \mathfrak s )$. 

\begin{enumerate}
\item 
\[
 -\kappa (Y , \mathfrak s )  \le 
 \begin{cases}
 b_2^+ (W) +\frac{\sigma (W)}{8}  -1 &\quad \text{if $b_2^+ (W)$ is odd,} \\
 b_2^+ (W) +\frac{\sigma (W)}{8} -2 &\quad \text{if $b_2^+ (W)$ is even and $b_2^+ (W) >0$ .} 
 \end{cases}
 \]
 
 \item 
 Furthermore if $(Y, \mathfrak s )$ is Floer $K_G$ split, then we have 
 \[
 \kappa (Y, \mathfrak s ) \le 
 \begin{cases}
 b_2^- (W) -\frac{\sigma (W)}{8} -1 &\quad \text{if $b_2^- (W)$ is odd ,} \\
 b_2^- (W) -\frac{\sigma (W)}{8} -2 &\quad \text{if $b_2^- (W)$ is even and $b_2^- (W) >0$. } 
 \end{cases}
 \]
 
 \end{enumerate}
 \end{corollary} 
 
 \begin{remark}
 Theorem \ref{estimate2} does not hold 
 for a cobordism $W$ with $b_2^+ (W) =0$ or $b_2^- (W) =0$. Likewise for a spin 4-manifold $(W , \mathfrak s_W )$ bounded by 
  $(Y, \mathfrak s )$ with $b_2^+ (W) =0$ (resp.~ 
 $b_2^- (W) =0$), we only have the inequality $-\kappa (Y, \mathfrak s ) \le b_2^+ (W) +\frac{\sigma (W)}{8} $ 
 (resp.~$\kappa (Y, \mathfrak s ) \le b_2^- (W) -\frac{\sigma (W)}{8}$) (as in the case where $W=S^3 \times [0,1]$ or $W =D^4$). 
 \end{remark}
 
The proofs of Theorems \ref{estimate1} and \ref{estimate2} are based on the finite dimensional approximation of the Seiberg-Witten map for  a cobordism $W$ with 
 certain boundary conditions developed in \cite{Mano1}. Let $(V_i )^{\tau}_{\mu}$  be the 
 finite dimensional vector spaces 
 used in the definitions of the approximated Seiberg-Witten flows for $(Y_i ,\mathfrak s_i )$, 
 and $I_i$ be the Conley indices $I_{\nu}$ used in the definitions of 
 $SWF (Y_i , \mathfrak s_i )$ $(i=0,1)$. Then we obtain a $G$ equivariant map  
$f : Z_0 \to Z_1$ where $Z_0$ and $Z_1$ are defined (up to stable equivalence) as follows (\cite{Mano4}).
 \begin{align*}
 &Z_0 =(\bh^{N +\ind_{\bh} \mathcal D_{W} (\mathfrak s_W ) -\dim (V_0 )^0_{-\nu} (\bh )})^+ \wedge 
 (\widetilde{\br}^{M-\dim_{\br} (V_0 )^0_{-\nu} (\widetilde{\br} ) } )^+ \wedge I_0 \\
 &Z_1 =(\bh^{N-\dim_{\bh} (V_1 )^0_{-\nu} (\bh )} )^+ \wedge (\widetilde{\br}^{M-\dim_{\br} (V_1 )^0_{-\nu} (\widetilde{\br} )
  +b_2^+ (W)} )^+ \wedge I_1 
 \end{align*}
 for some $M$, $N$ (we may assume that $M$ is even).
Here $Z_0$ and $Z_1$ are spaces of type SWF at level $M$ and $M+b_2^+ (W)$ respectively, and 
the restriction of $f$ to the $S^1$ fixed point sets induces the map 
\[
f^{S^1} : Z_0^{S^1} \cong S^M  \to  Z_1^{S^1} \cong (\widetilde{\br}^{b_2^+ (W)})^+ \wedge S^M
\]
which is a $G$ homotopy equivalence if $b_2^+ (W) =0$. Furthermore the restriction of $f$ to the $G$ fixed point sets induces a $G$ homotopy equivalence
$f^G : Z_0^G \to Z_1^G$. Hence if $b_2^+ (W)$ is even, by Lemma \ref{estimate for SWF} we obtain 
\[
k(Z_0 ) \le k(Z_1 ) +\frac 12 b_2^+ (W) .
\]

If $I_i$ is at odd level, then $I_i$ and $M$ in the description of $Z_i$ are replaced with $\Sigma^{\widetilde{\br}} I_i$ and $M-1$ respectively. Then we have 
\begin{align*}
&k(Z_0 ) =\ind_{\bh} \mathcal D_W (\mathfrak s_W ) -\dim_{\bh} (V_0 )^0_{-\nu} (\bh ) +k(I_0 ')  +N
\\
&k(Z_1 ) =-\dim_{\bh} (V_1)^0_{-\nu} (\bh ) +k(I_1 ') +N
\end{align*}
where $I_i ' =I_i$ (resp.~$\Sigma^{\widetilde{\br}} I_i $) 
if $I_i$ is at even level (resp.~odd level), since $k(\Sigma^{\widetilde{\bc}} X ) =k(X)$ and $k(\Sigma^{\bh} X ) =k(X) +1$ for any space $X$ of type SWF 
at even level 
(\cite{Mano4}). 
 It follows from the above inequality and the definition of $\kappa (Y_i, \mathfrak s_i )$ that 
 \begin{align*}
 \kappa (Y_0 , \mathfrak s_0 ) -\kappa (Y_1 , \mathfrak s_1 ) &\le b_2^+ (W) -\ind_{\bc} \mathcal D_W (\mathfrak s_W )  + 
 n(Y_1 , \mathfrak s_1 , g_1 ) -n(Y_0 , \mathfrak s_0 , g_0 )  \\
 &= b_2^+ (W) +\frac{\sigma (W) }{8} 
 \end{align*}
since $n (Y_1 , \mathfrak s_1 , g_1 ) -n (Y_0 , \mathfrak s_0 , g_0 ) =\ind_{\bc} \mathcal D_W (\mathfrak s_W ) +\frac{\sigma (W)}{8}$ 
(Here we assume that $\ker \mathcal D_{Y_i} (\mathfrak s_i ) =0$ by 
 using the perturbed Chern-Simons-Dirac functional if necessary, which ensures that 
$n (-Y_i , \mathfrak s_i , g_i ) =-n(Y_i , \mathfrak s_i , g_i )$.)
If $b_2^+ (W)$ is odd, by using $W\sharp S^2 \times S^2$ instead of $W$ we deduce 
\[
\kappa (Y_0  , \mathfrak s_0 ) -\kappa (Y_1 , \mathfrak s_1 ) \le b_2^+ (W) +\frac{\sigma (W)}{8} +1 .
\]
We obtain the other inequalities in Theorem \ref{estimate1} by applying the above argument to the cobordism $(-W, \mathfrak s_W )$ from 
$(Y_1 , \mathfrak s_1)$ to $(Y_0 , \mathfrak s_0 )$. 
 If $(Y_0, \mathfrak s_0 )$ is Floer $K_G$ split and $b_2^+ (W)$ is even and nonzero, then Lemma \ref{estimate for SWF} shows that   
$ k(Z_0 ) + 1 \le k(Z_1 ) +\frac{b_2^+ (W)}{2}$, from which we deduce 
\begin{align*}
\kappa (Y_0 , \mathfrak s_0 ) -\kappa (Y_1 , \mathfrak s_1 )  &\le 
b_2^+ (W) -2  -\ind_{\bc} \mathcal D_W (\mathfrak s_W ) \\
 &+n(Y_1 , \mathfrak s_1 , g_1 ) -n(Y_0 , \mathfrak s_0 , g_0 )  
= b_2^+ (W) +\frac{\sigma (W)}{8} -2 .
\end{align*}
We obtain the other inequalities in Theorem \ref{estimate2} by replacing $W$ with $W\sharp S^2 \times S^2$ or $-W$ as before. 

\section{Proof of Theorem \ref{main2}}

To prove Theorem \ref{main2} we first consider the relation between the $\kappa$ and the $\overline{\mu}$ invariants of a spherical 3-manifold. 

\begin{proposition}\label{PSC} \cite{Mano4}
Let $(Y, \mathfrak s)$ be a rational homology 3-sphere with $\spin$ structure that admits a metric $g$ of positive scalar curvature. Then 
$\kappa (Y, \mathfrak s ) = -n( Y, \mathfrak s , g )$ and $(Y, \mathfrak s )$ is Floer $K_G$ split. 
\end{proposition}

In the above case, the critical point set of the approximated Seiberg-Witten flow on $V$ only consists of the unique reducible element. Hence 
we have $I_{\nu} = (V^0_{-\nu} )^+$  and $\kappa (Y, \mathfrak s ) =k(S^0 ) -n(Y, \mathfrak s , g ) =-n(Y, \mathfrak s, g )$. 
Furthermore the ideal $\mathfrak J (I_{\nu} )$ (or $\mathfrak J (\Sigma^{\widetilde{\br}} I_{\nu} )$ if $I_{\nu}$ is at odd level) is generated by 
$z^k$ for some $k$, since $\mathfrak J (\Sigma^{\widetilde{\bc}} X)=\mathfrak J (X)$ and $\mathfrak J (\Sigma^{\bh} X ) =z\mathfrak J (X )$ for a 
space $X$ of type SWF at even level (\cite{Mano4}). Hence $(Y, \mathfrak s )$ is Floer $K_G$ split. 

For a rational homology 3-sphere $(Y, \mathfrak s )$ with metric $g$, let $(W, \mathfrak s_W )$ be a compact $\spin$ 4-manifold bounded by 
$(Y, \mathfrak s )$ and $g_W$ be a metric of $W$ extending $g$. Then the index of the Dirac operator
$\mathcal D_W (\mathfrak s_W )$ of $(W, \mathfrak s_W )$ associated with $g_W$ and the 
signature of $W$ are described as 
\begin{align*}
\ind_{\bc} \mathcal D_W (\mathfrak s_W ) &=-\frac{1}{24} \int_W p_1  -\frac 12 (\dim \ker \mathcal D_Y (s ) +\eta^{\mathrm{Dir}} (Y, s , g ) ) ,  \\
\sigma (W) &=\frac 13 \int_W p_1 -\eta^{\mathrm{sign}} (Y, g ) .
\end{align*} 
 It follows that 
 \[ (1) \quad 
 n (Y, \mathfrak s , g ) =-\frac 18 (4\dim \ker \mathcal D_Y (\mathfrak s )  +4\eta^{\mathrm{Dir}} (Y, \mathfrak s , g )  +
 \eta^{\mathrm{sign}}  (Y, g ) ). 
 \]
  If $g$ is a metric of positive scalar curvature, $\dim \ker \mathcal D_Y (\mathfrak s ) =0$.  
 For a spherical 3-manifold with $\spin$ structure $(S, \mathfrak s )$ equipped with the standard metric $g$, 
  we have the following equation as a special case of Proposition \ref{U2}, which is deduced from Proposition \ref{orbifold index} and Proposition \ref{U3}. 
  \[ (2) \quad 
  \overline{\mu} ( S, \mathfrak s ) =w (S, cS , \mathfrak s_{cS} ) =-\frac 18 (4\eta^{\mathrm{Dir}} (S, \mathfrak s , g ) +\eta^{\mathrm{sign}} (S, g )) 
  \]
  where $\mathfrak s_{cS}$ is a $\spin$ structure on the cone $cS$ over $S$ extending $\mathfrak s $.  It follows from (1), (2) and Proposition \ref{PSC} that 
  $\kappa (S, \mathfrak s ) =-\overline{\mu} (S, \mathfrak s )$. 
  
\begin{remark}\label{sign}
To describe the above equation we assume that the Clifford multiplication $c$ of the volume form $\mathrm{vol}_Y$ of $Y$ satisfies $c (\mathrm{vol}_Y ) =-1$ 
as in \cite{Ni}. If we define $c$ so that it satisfies $c(\mathrm{vol}_Y ) =1$ as in  \cite{Mano1}, \cite{Mano3}, \cite{Mano4}, 
the sign of $\eta^{\mathrm{Dir}} (Y, \mathfrak s, g )$ 
should be changed and the first equation in Proposition \ref{U3} should be replaced with 
$\eta^{\mathrm{Dir} } (S, \mathfrak s , g ) =-2\delta^{\mathrm{Dir}} (S, \mathfrak s )$. 
\end{remark}

 Suppose that a rational homology 3-sphere with $\spin$ structure $(Y, \mathfrak s )$ bounds a compact $\spin$ 4-orbifold $(X, \mathfrak s_X )$ 
 with $\mathfrak s_X |_Y =\mathfrak s$. Let 
 $(S_i , \mathfrak s_i )$ $(i=1, \dots , n )$ be spherical 3-manifolds, which are the links of all isolated singularities of $X$, where a $\spin$ structure $\mathfrak s_i$ is induced from $\mathfrak s_X$. 
Putting $X_0 =X \setminus \cup_{i=1}^n \Int cS_i$, we have a compact $\spin$ 4-manifold with 
$\partial (X_0 , \mathfrak s_{X_0} ) =(Y , \mathfrak s ) \cup \cup_{i=1}^n (-S_i , \mathfrak s_i )$ whose $\spin$ structure $\mathfrak s_{X_0}$ is 
induced from $\mathfrak s_X$. 
Furthermore $(X_0 , \mathfrak s_{X_0} )$ contains a $\spin$ submanifold $W$
of codimension $0$ with 
$\partial W = \sharp_{i=1}^n S_i \cup \cup_{i=1}^n (-S_i )$ which is constructed from 
the collars of $S_i$ by attaching 1-handles in $X_0$. 
Let $X_0 ' = X_0 \setminus \Int W$ and $\mathfrak s_{X_0 '}$ be a $\spin$ structure induced from $\mathfrak s_X$. Then 
$(X_0 ' , \mathfrak s_{X_0 '} )$ is a $\spin$ cobordism from $(S_0, \mathfrak s_0 )
 = (\sharp_{i=1}^n S_i , \sharp_{i=1}^n \mathfrak s_i )$ to $(Y, \mathfrak s )$. 
We note that $S_0$ admits a metric $g_0$ of positive scalar curvature (\cite{GL2}, \cite{SY}). 
Thus if we choose a metric $g_W$ on $W$ extending $g_0$ and the standard metrics $g_i$ of $S_i$ $(i\ge 1)$, 
we obtain a $G$ equivariant map $f: Z_0 \to Z_1$ for the cobordism $W$ from $\cup_{i=1}^n S_i$ to $S_0$ 
of the following form (up to stabilization) by a procedure similar to those in \cite{Mano3}, \cite{Mano4} since the Conley index for a disjoint union of $S_i$ 
is the smash product of the Conley indices of $S_i$ (\cite{Mano3}). 

\begin{align*}
&Z_0 =(\bh^{N +\ind_{\bh} \mathcal D_W (\mathfrak s_W ) -\sum_{i=1}^n \dim_{\bh} (V_i )^0_{-\nu} (\bh )})^+ 
\wedge (\widetilde{\br}^{M-\sum_{i=1}^n \dim_{\br} (V_i )^0_{-\nu}  (\br )} )^+ \wedge (\wedge_{i=1}^n I_i ) \\
&Z_1 =
(\bh^{N-\dim_{\bh} (V_0 )^0_{-\nu} (\bh )} )^+  \wedge (\widetilde{\br}^{M-\dim_{\br} (V_0 )^0_{-\nu}  (\widetilde{\br} ) +b_2^+ (W) } )^+ \wedge I_0 
\end{align*}
where $(V_i )^0_{-\nu}$ and $I_i$ are the vector space and the Conley index for $S_i$ defined as above $(0\le i \le n )$. 
Since $I_i =(V_i )^0_{-\nu}$ and $b_2^+ (W) =0$, 
$Z_0$ and $Z_1$ are of type SWF at level $M$ (which can be assumed to be even), 
$f^{S^1} : Z_0^{S^1} \to Z_1^{S^1}$ is a $G$ homotopy equivalence. Furthermore 
$k(Z_0 ) =N +\ind_{\bh} \mathcal D_W (\mathfrak s_W )$ and $k(Z_1 ) =N$. Hence 
applying Lemma \ref{estimate for SWF} to $f$ we deduce that $\ind_{\bh} \mathcal D_W (\mathfrak s_W ) \le 0$. 
We note that since $\sigma (W) =0$,  we have 
\begin{align*}
&\ind_{\bc} \mathcal D_W (\mathfrak s_W ) =-\frac{1}{24} \int_W p_1 -\frac 12 (\eta^{\mathrm{Dir}} (S_0 , \mathrm s_0 , g_0 ) 
+\sum_{i=1}^n \eta^{\mathrm{Dir}} (-S_i , \mathfrak s_i , g_i )  )  
\\
&=
-\frac 18 (4\eta^{\mathrm{Dir}} (S_0, \mathfrak s_0 , g_0 ) +\eta^{\mathrm{sign}} (S_0, g_0 )  + 
\sum_{i=1}^n (\eta^{\mathrm{Dir}} (-S_i , \mathfrak s_i , g_i ) +\eta^{\mathrm{sign}} (-S_i , g_i ) ) ) \\
&= n(S_0, \mathfrak s_0 , g_0 ) -\sum_{i=1}^n  n(S_i , \mathfrak s_i , g_i ) 
\end{align*}
By a similar map defined for a cobordism $-W$ from $S_0$ to $\cup_{i=1}^n S_i$, we have 
$\ind_{\bh} \mathcal  D_{(-W)} (\mathfrak s_W )=-\ind_{\bh} \mathcal D_W (\mathfrak s_W ) \le 0$ since $b_2^+ (-W) =b_2^- (W)=0$. 
 It follows that $\ind_{\bc} \mathcal D_W (\mathfrak s_W ) =2\ind_{\bh} \mathcal D_W (\mathfrak s_W ) =0$ and 
\[ (3) \quad 
\kappa (S_0 , \mathfrak s_0 ) =-n(S_0 , \mathfrak s_0 , g_0 ) =-\sum_{i=1}^n n(S_i , \mathfrak s_i , g_i ) =-\sum_{i=1}^n \overline{\mu} (S_i , \mathfrak s_i ) .
\]

Hereafter we assume that $(Y, \mathfrak s )$ is a Seifert rational homology 3-sphere with $\spin$ structure and choose  
 the $\spin$ 4-orbifolds $(X^{\pm} , \mathfrak s_{X^{\pm}} )$
bounded by $(Y, \mathfrak s )$ satisfying the conditions in Proposition \ref{U2}. Let $(S_i^{\pm} , \mathfrak s_i^{\pm} )$ 
$( i=1 ,\dots , n^{\pm})$ 
be the spherical 3-manifolds which are the links of the isolated singularities of $X^{\pm}$, where $\mathfrak s_i^{\pm}$ is induced from 
$\mathfrak s_{X^{\pm}}$. 
Subtracting all the cones $cS_i^{\pm}$ and the $\spin$ 4-manifold of codimension $0$, which is a $\spin$ cobordism from 
$\cup_{i=1}^{n^{\pm}} (S_i^{\pm} , \mathfrak s_i^{\pm} )$ to $(S_0^{\pm} , \mathfrak s_0^{\pm} ) =
(\sharp_{i=1}^{n^{\pm}} S_i^{\pm} , \sharp_{i=1}^{n^{\pm}} \mathfrak s_i^{\pm} )$ constructed as above from $X^{\pm}$, 
we obtain compact $\spin$ cobordisms $(X_0^{\pm} , \mathfrak s_{X_0^{\pm}} )$ from 
$(S_0^{\pm} , \mathfrak s_0^{\pm} )$ to $(Y, \mathfrak s )$ satisfying 
\[
b_1 (X_0^{\pm} ) =b_1 (X^{\pm} ) =0, \ \  b_2^+ (X_0^+ ) =b_2^+ (X^+ ) \le 1, \ \  b_2^- (X_0^-   ) =b_2^- (X^- ) \le 1 . 
\]

Since $(S_0^{\pm} ,  \mathfrak s_0^{\pm} )$ is Floer $K_G$ split, we have the following inequalities by applying Theorem \ref{estimate1} and 
Theorem \ref{estimate2}. 

\begin{align*}
(4) \quad  
-b_2^-  (X_0^{\pm} ) +\frac{\sigma (X_0^{\pm} ) }{8} -1 &\le \kappa  (S_0^{\pm} , \mathfrak s_0^{\pm}  ) -\kappa (Y, \mathfrak s ) \\
\le b_2^+ (X_0^{\pm} ) +\frac{\sigma (X_0^{\pm}  ) }{8} -1 
\end{align*}
unless $b_2^+ (X_0^{\pm}) =0$.  If $ b_2^+ (X_0^+) =0$, we replace $X_0^+$ with 
$X_0^+ \sharp S^2 \times S^2$ to obtain the inequality 
$\kappa (S^+  , \mathfrak s_{S^+} ) -\kappa (Y, \mathfrak s ) \le \sigma (X_0^+ )/8$ by applying Theorem \ref{estimate2}. 
If $b_2^- (X_0^- )$ is even (in fact $0$), then the first inequality of (4) for $S^-$ can be replaced with 
\[
(5) \quad 
-b_2^- (X_0^-) +\frac{\sigma (X_0^- )}{8}  \le \kappa (S^- , \mathfrak s_{S^- } ) -\kappa (Y,\mathfrak s ) . 
\]

On the other hand it follows from Propositions \ref{orbifold index}, \ref{U3}, \ref{U2} and the equations (1), (2), (3) that 
\begin{align*}
(6) \quad 
&\overline{\mu} (Y, \mathfrak s ) = w(Y, X^{\pm} , \mathfrak s_{X^{\pm}} ) \\
&=\frac 18 (\sigma (X_0^{\pm} ) - \sum_{i=1}^n (
4\eta^{\mathrm{Dir}} (S_i^{\pm} , \mathfrak s_i^{\pm} , g_i^{\pm}   ) +\eta^{\mathrm{sign}} (S_i^{\pm} , g_i^{\pm} ) )) \\
&=\frac 18 \sigma (X_0^{\pm} ) +\sum_{i=1}^n n(S_i^{\pm} , \mathfrak s_i^{\pm} , g_i^{\pm} ) =\frac 18 \sigma (X_0^{\pm} ) -\kappa 
(S_0^{\pm}, \mathfrak s_0^{\pm}  )
\end{align*}

Therefore by subtracting (6) from (4), we obtain

\[
 (7) \quad 
-2 \le -b_2^- (X_0^- ) -1 \le -\kappa (Y, \mathfrak s ) -\overline{\mu} (Y, \mathfrak s ) \le b_2^+ (X_0^+ ) -1 \le  0 
\]
if $b_2^+ (X_0^+ ) =1$. We replace $X_0^+$ with 
$X_0^+  \sharp S^2 \times S^2$ or 
$X_0^+ \sharp 2S^2 \times S^2$ if $b_2^+ (X_0^ + ) =0$ to obtain the last inequality by
applying Theorem \ref{estimate2}. 
If $Y$ is a Seifert integral homology 3-sphere (with the unique $\spin$ structure), 
$\kappa (Y) +\overline{\mu} (Y) =0$ or $2$ since both $\kappa (Y)$ and $\overline{\mu} (Y)$ are integers and 
the integral lifts of the Rokhlin invariant. 

Suppose that a Seifert rational homology 3-sphere $(Y, \mathfrak s )$ is Floer $K_G$ split. Then by Theorem \ref{estimate2} we obtain 
\[
0 \le -b_2^-  (X_0^- ) +1 \le -\kappa (Y, \mathfrak s ) -\overline{\mu} (Y , \mathfrak s) 
\]
if $b_2^- (X_0^- ) =1$ in addition to (7). 
We obtain the same estimate by (5) or applying Theorem \ref{estimate2} to $X_0^- \sharp S^2 \times S^2$ if 
$b_2^- (X_0^- ) =0$.  Thus we have $\kappa (Y , \mathfrak s ) + \overline{\mu} (Y , \mathfrak s ) =0$. 

For a general rational homology 3-sphere with $\spin$ structure $(Y, \mathfrak s )$,  
 it is proved that $\kappa (Y, \mathfrak s ) +\kappa (-Y , \mathfrak s ) \ge 0$ (\cite{Mano4}). Note that 
 $\kappa (Y, \mathfrak s)$ is not necessarily the same as $-\kappa (-Y, \mathfrak s )$, while
we have $\overline{\mu} (-Y, \mathfrak s ) =
-\overline{\mu} (Y, \mathfrak s )$.  
If $Y$ is a Seifert rational homology 3-sphere, the above result shows that 
$0\le \kappa (Y, \mathfrak s ) +\overline{\mu} (Y, \mathfrak s ) \le 2$, 
$0 \le \kappa (-Y, \mathfrak s ) -\overline{\mu} (Y, \mathfrak s ) \le 2$, 
and hence $0\le \kappa (Y, \mathfrak s ) + \kappa (-Y, \mathfrak s ) \le 4$. If $Y$ is a Seifert integral homology 3-sphere, we also have 
$\kappa (-Y) +\overline{\mu} (-Y) =\kappa (-Y ) -\overline{\mu} (Y) =0$ or $2$. 
It follows that $\kappa (Y) + \kappa (-Y) =0$,$2$, or $4$.
If $Y$ has a fibration such that one of the multiplicities is even and has a positive degree, then 
 we apply the above inequalities for $Y$ to $X^- =X^+$ satisfying $b_2^-  (X^{\pm} ) =0$ and $b_2^+ (X^{\pm} ) =1$, which is chosen in Proposition \ref{U2}. 
 It follows from (5), (6), (7) that $\kappa (Y ,\mathfrak s) =-\overline{\mu}  (Y, \mathfrak s)$.

This completes the proof of Theorem \ref{main2}. 

\begin{example}
We give a list of the values of $\kappa$, $\overline{\mu}$
together with $\beta$ and $\underline d$ of some Brieskorn homology 3-spheres. 
The computations of the $\kappa$ invariants below are due to Manolescu \cite{Mano4}.  
It is pointed out in \cite{Mano4} that $\pm \Sigma (2,3,12n+1)$ and $\pm \Sigma (2,3,12n+5 )$ are Floer $K_G$ split, while 
$\pm \Sigma (2,3,12n-1)$ and $\pm \Sigma (2,3,12n-5 )$ are not.

 \qquad 
\begin{tabular}{ccccc}
 & \quad $\kappa$ \quad & \quad $\overline{\mu}$  \quad & \quad $\beta$ \quad & \quad $\underline d$ \quad \\
 \hline \hline 
 $\ \Sigma (2,3, 12n-1 )$ & \quad 2 \quad & \quad 0  \quad & \quad 0 \quad & \quad 0 \quad \\
 \hline
 $-\Sigma (2,3,12n-1 )$ & \quad 0 \quad & \quad 0 \quad & \quad 0 \quad & \quad 0 \quad \\
 \hline 
 $\ \Sigma (2,3,12n+1)$ & \quad 0 \quad & \quad 0 \quad & \quad 0 \quad & \quad 0 \quad \\
 \hline
 $-\Sigma (2,3,12n+1)$ & \quad 0 \quad & \quad 0 \quad &\quad  0 \quad & \quad 0 \quad \\
 \hline 
 $\ \Sigma (2,3,12n-5)$ & \quad 1 \quad & \quad 1  \quad & \quad -1 \quad &\quad  -2 \quad \\
 \hline 
 $-\Sigma (2,3,12n-5 )$ & \quad 1 \quad &\quad  -1  \quad & \quad 1 \quad & \quad 2 \quad \\
 \hline 
 $\ \Sigma (2,3,12n+5)$ &\quad  1 \quad & \quad -1 \quad  & \quad 1 \quad & \quad 2 \quad \\
 \hline 
 $-\Sigma (2,3,12n+5)$ & \quad -1\quad & \quad 1 \quad & \quad -1 \quad &\quad  -2  \quad \\
 \hline 
 \end{tabular}
 
 \end{example}

\begin{remark}\label{addenda}
If $Y$ is a Seifert fibration of the form $\{b, (a_1 , b_1) , \dots  , (a_n , b_n )\}$, then $Y$ is represented by a framed link $\cup_{i=0}^n K_i$, where 
$K_0$ is an unknot with framing $b$ and $K_i$ $(i\ge 1 )$ is a meridian of $K_0$ with framing $-b_i /a_i$. 
Let $h$ be a meridian of $K_0$ (a general fiber of $Y$) and $g_i$ be a meridian of $K_i$ $(i\ge 1 )$. Then 
a $\spin$ structure $\mathfrak s$ on $Y$ 
corresponds to a homomorphism $c: H_1 (S^3 \setminus \cup_{i=0}^n K_i , \bz ) \to \bz_2$ satisfying 
\begin{align*}
& a_i c(g_i ) + b_i c(h) \equiv a_i b_i \pmod 2, \\
& \sum_{i=1}^n c(g_i ) +bc(h) \equiv b \pmod 2 
\end{align*}
(\cite{U2}. The description of the Seifert invariants in \cite{U2} is different from that of this paper, but the above condition is still valid.)
Suppose that all $a_i$ are odd and $\sum_{i=1}^n b_i \equiv b \pmod 2$ (such a case only occurs when $|H_1 (Y, \bz )|$ is even). 
Then we have a $\spin$ structure $\mathfrak s$ on $Y$ corresponding to $c$ that satisfies $c(h) \equiv 0$, $c(g_i ) \equiv b_i$ $(i\ge 1 )$. 
 In this case $X^\pm$ in Proposition \ref{U2} can be also chosen so that 
 $X^+ =X^-$ satisfying $b_2^+ (X^{\pm} ) =1$ and $b_2^- (X^{\pm} ) =0$ if $\deg Y >0$, and 
 $b_2^- (X^{\pm} ) =1$ and $b_2^+ (X^{\pm} ) =0$ if $\deg Y <0$ (\cite{U2}). 
 Thus it follows from the proof of Theorem 2 that $\kappa (Y, \mathfrak s ) =-\overline{\mu} (Y, \mathfrak s )$ if $\deg Y >0$. 
 \end{remark}

\begin{remark}
 If a Seifert rational homology 3-sphere with $\spin$ structure $(Y, \mathfrak s )$ bounds a negative definite $\spin$ 4-manifold $W$, the inequalities in Theorem \ref{main1} coincide with those in \cite{U3}. 
Furthermore 
 for a Seifert integral homology 3-sphere, the second inequality in Theorem \ref{main1} 
 coincides with those in terms of the $\beta$ invariant (\cite{Mano3}), and in terms of the $\underline{d}$ invariant 
(\cite{HM}) under the relations 
$\beta (Y) =-\overline{\mu} (Y)$ (\cite{Dai}, \cite{Stoff}) and $\underline{d} (Y) =-2\overline{\mu} (Y)$ 
 (\cite{DM}).
\end{remark}

\begin{remark}
We have $\kappa (Y) +\overline{\mu} (Y) =0$ or $2$ for a Seifert integral homology 3-sphere $Y$
by Theorem \ref{main2}. 
If $\kappa (Y) =-\overline{\mu} (Y)$, then Corollary \ref{cor} shows that 
$\overline{\mu} (Y ) \le b_2^+ (W) +\frac{\sigma (W)}{8} -1$ for a compact $\spin$ 4-manifold $W$ bounded by $Y$, 
which give the better estimate than
Theorem \ref{main1}. If $\kappa (Y) = 2-\overline{\mu} (Y)$, the estimate in Theorem \ref{main1} is slightly better than or the same as 
that of $\overline{\mu} (Y)$ given in Corollary \ref{cor}. 
\end{remark}

J.Lin \cite{JL} proves other  constraints on the intersection forms of $\spin$ 4-manifolds bounded by homology 3-spheres in terms of $\mathrm{KO}$ invariants, 
which give estimates better than those in Theorem \ref{main1} and Theorems \ref{estimate1}, \ref{estimate2} in some cases. 
In \cite{JL} several estimates for $\spin$ 4-manifolds bounded by certain Brieskorn homology 3-spheres are also given by considering the 
closed $\spin$ 4-orbifolds as in \cite{FF}. We obtain similar estimates by using Proposition \ref{U2} as follows. 
Suppose that a Seifert rational homology 3-sphere $Y$ (with $\spin$ structure) bounds a $\spin$ 4-manifold $W$. 
Let $Z = X^- \cup (-W)$ be a closed $\spin$ 4-orbifold, where $X^-$ is a $\spin$ 4-orbifold bounded by $Y$ chosen in Proposition \ref{U2}. 
Then we have a $G$ equivariant map 
\[
f: (\bh^{\ind_{\bh} \mathcal D_{-Z} } )^+ \to (\widetilde{\br}^{b_2^+ (-Z )} )^+
\]
whose restriction to the $G$-fixed point set is a homotopy equivalence (we choose $Z$ according to our sign convention). 
It is proved in \cite{JL} that 
there is no such map of the form $f: (\bh^{4\ell} )^+ \to (\widetilde{\br}^{8\ell +2} )^+$ if $\ell >0$. Note that 
$b^- (X^- ) \le 1$ and $\ind_{\bc} \mathcal D_{-Z} =\overline{\mu} (Y) -\frac{\sigma (W)}{8}$. It follows that
if $\overline{\mu} (Y) - \frac{\sigma (W)}{8} >0$ and divisible by 8, then 
\[
\overline{\mu} (Y) \le b_2^+ (W) +\frac{\sigma (W)}{8} -2 ,
\]
which gives a better estimate than those by Theorem \ref{main1} or Corollary \ref{cor}.

 \end{document}